\documentclass[journal]{IEEEtran}
\usepackage{amsmath}
\usepackage{amsmath,amsfonts}
\usepackage{cite}
\usepackage{graphicx,amssymb}
\usepackage{amsmath,amsfonts,amssymb}

\usepackage[usenames]{color}
\usepackage{algorithm}
\usepackage{algorithmic}
\usepackage{epstopdf}
\usepackage{multirow}
\usepackage{stfloats}
\usepackage{subfigure}
\usepackage{url}

\newcommand{\Tr}{\operatorname{Tr}}

\begin{document}
 
\title{Testing Link Fidelity in a Quantum Network using Operational Form of Trace Distance with Error Bounds}

\author{J.T. M. Campbell~\IEEEmembership{Student Member,~IEEE}, N. Marchetti,~\IEEEmembership{Senior Member,~IEEE}, John Dooley,~\IEEEmembership{Member,~IEEE}, and I. Dey,~\IEEEmembership{Senior Member,~IEEE}

\thanks{J. T. M. Campbell is with Department of Electronic Engineering, National University of Ireland, Maynooth (E-mail: john.campbell.2023@mumail.ie).}
\thanks{N. Marchetti is with Department of Electrical and Electronic Engineering, Trinity College Dublin, Ireland (E-mail: nicola.marchetti@tcd.ie).}
\thanks{John Dooley is with Department of Electronic Engineering, National University of Ireland, Maynooth (Email: John.Dooley@mu.ie)}
\thanks{I. Dey is with Walton Institute for Information and Communications Systems Sciences, Waterford, Ireland (Email: indrakshi.dey@waltoninstitute.ie)}
\thanks{This material is based upon work supported by Science Foundation Ireland (SFI) and is co-funded under the European Regional Development Fund under Grant Number 13/RC/2077 and 13/RC/2077-P2.}
}

%\markboth{Journal of \LaTeX\ Class Files,~Vol.~6, No.~1, January~2007}%
%{Shell \MakeLowercase{\textit{et al.}}: Bare Demo of IEEEtran.cls for Journals}

\maketitle

\begin{abstract}
Quantum state comparison, utilizing metrics like fidelity and trace distance, underpins the assessment of quantum networks within quantum information theory. While recent research has expanded theoretical understanding, incorporating error analysis and scalability considerations remains crucial for practical applications. The primary contribution of this letter is to address these gaps by deriving the novel operational trace distance for multi-node networks, establishing a trace distance vs. fidelity benchmark incorporating error bounds, and bridging quantum operations with tensor network analysis.  We further explore the application of tensor network tools to quantum networks, offering new analytical avenues.  This comprehensive approach provides a robust framework for evaluating quantum network performance under realistic error conditions, facilitating the development of reliable quantum technologies.

\end{abstract}

\IEEEpeerreviewmaketitle
\vspace*{-7mm}

\section{Introduction}
\label{S1}

Quantum state comparison plays a fundamental role in quantum information theory, enabling the assessment of similarity or dissimilarity between quantum states. Various measures, such as fidelity and trace distance, have been proposed to quantify the degree of similarity between states. Fidelity, defined as the overlap between two states, captures their similarity, while trace distance quantifies the maximum distinguishability. These measures serve as benchmarks for evaluating the effectiveness and performance of quantum networks, which are crucial for the implementation of quantum communication and computation protocols.

In recent years, there has been significant progress in understanding the applications and implications of fidelity and trace distance as benchmarks for state comparison. The literature has explored their interpretation, significance in preserving quantum coherence and entanglement, and in relation to various mathematical properties such as logarithmic negativity and entropy \cite{Helsen2023, Ruskai2001, Schroeder2022, Bennett1992, Nielsen2011, Berta2019, Wilde2013, Tomamichel2013}. However, there are still notable gaps in the bibliography, particularly regarding the incorporation of errors into quantum network benchmarks and the scalability of these benchmarks for example in review of the literature \cite{Helsen2023, Berta2019, Wilde2013, Tomamichel2013}. While works such as \cite{Helsen2023} exist that present a benchmarking procedure for quantum networks, including such bounds would enhance the assessment of the performance and reliability of the benchmarked networks.

While [2] provides a strong foundation in the mathematics of quantum information, the inclusion of upper and lower error bounds would further solidify the analytical framework for quantum information processing. Similarly, [3], [4], [5] and [7], offer valuable insights into quantum communication and computation but would benefit from an explicit treatment of error bounds. Incorporating error bounds would enhance the assessment of protocol robustness, such as in the case of Clauser-Horne-Shimony-Holt (CHSH)-based protocols.

In the realm of fidelity, [6] and [8] investigate properties of fidelity; however, a comprehensive understanding of fidelity of recovery necessitates the consideration of error bounds. This extends to the analysis of continuous variable quantum networks in [9], where error bound analysis would provide deeper insights and connections to classical complex network error models. The practical implementation of quantum networks highlights their inherent susceptibility to errors. A rigorous analysis, incorporating error bounds into fidelity and trace distance benchmarks, is essential for evaluating state similarity and reliable information processing in real-world quantum networks. Addressing the current gap in the literature is crucial for understanding the limitations and ensuring the robustness of these benchmarks under realistic conditions.

To bridge existing research gaps, this letter offers the following key contributions: a) We derive and implement a practical form of trace distance applicable to 2, 3, and multinode quantum networks, enabling rigorous assessment of state similarity. b) We investigate the impact of nearest and farthest node pairings on overall network fidelity, offering insights for optimizing network architecture. c) We establish a new benchmark based on the interpretation of the trace distance against fidelity relationship. Incorporating error bounds (from Fuchs and van de Graaf) enhances the benchmark's robustness for real-world applications. d) We also connect core quantum operations to tensor network diagrams, opening avenues for applying tensor network techniques to quantum information analysis. e) We adapt tools from tensor network analysis (e.g., Singular Value Decomposition, eigenvalue analysis) for the targeted study of quantum networks, and finally, we draw connections between these adapted tools and their counterparts in classical graph theory, providing a richer theoretical framework for understanding quantum network behavior.

\section{Derivation of an Operational Form of the Trace Distance}\label{S2}

Here we implement the use of the operational trace distance as a test of link fidelity within a simple quantum network consisting of two quantum nodes. We rely on the concept that each quantum node within a quantum network can be represented by its density matrix [2]; so essentially we are leveraging operational trace distance between density matrices to measure link fidelity. We start by defining a quantum channel $S : N, M$, where $N$ and $M$ refer to the Hilbert space dimensions of the input and output systems, respectively. For example, if we have a quantum channel $S: N, M$ where $N = 3$ and $M = 2$, it means that the input system has a Hilbert space dimension of 3, and the output system has a Hilbert space dimension of 2. 

We also introduce a measurement operator $\{\Lambda_x\}=\Lambda, \sum_x \Lambda_x = I$, where $x$ is a classical variable bit and $I$ represents the Identity operator. We consider that the measurement operator operates on the density matrices  $\rho \succeq 0$ and $\sigma \succeq 0$, where $\rho$ and $\sigma$ are the density matrices used for representing the two quantum nodes within the quantum network. If $\rho^2=\rho$, the density matrix represents a noiseless scenario, i.e. the quantum state of the involved quantum node exists on the surface of the Poincare / Bloch sphere [3]. It is worth-mentioning here that in quantum mechanics, a quantum state (corresponding to an energy level in a harmonic oscillator) is degenerate if it corresponds to two or more different measurable states of a quantum system. When degenerate or mixed, the density matrix will then exist within the Poincare/Bloch sphere. Furthermore, if $\rho$ and $\sigma$ are isolated in terms of quantum nodes, they will be in pure quantum states. If $\rho$ and $\sigma$ are connected through an interacting local channel, like in a quantum network, they will emerge as a Hartree-Fock wavefunction, $\psi$ of $\rho-\sigma$ mixed states.

\begin{figure}
    \centering
    \includegraphics[width=0.50\textwidth]{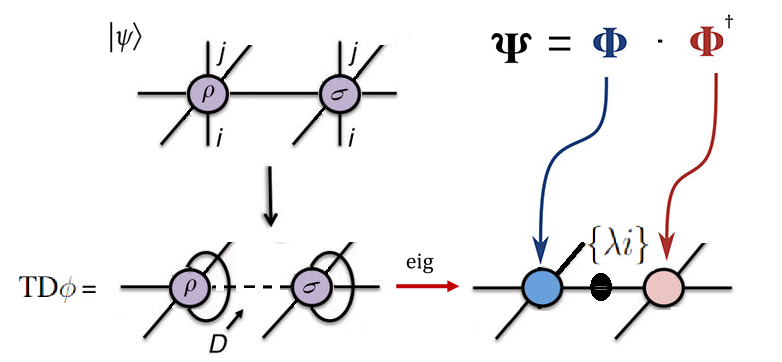}
    \caption{Tensor network diagrams of trace distance and trace norm of density matrices $\rho$ and $\sigma$, which represent a recoverable mixed-state of the wavefunction $\psi$ in a Hartree-Fock model}
    \label{model}
\end{figure}

\begin{figure}
    \centering
\includegraphics[width=0.45\textwidth]{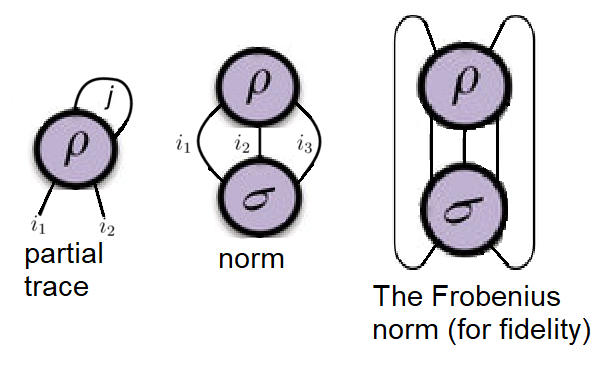}\\
%\centering
\vspace{-3mm}
\caption{Tensor Network Diagrams for Matrix Operations}
\vspace{-5mm}
\label{model1}
\end{figure}

In analyzing quantum networks, tensor network notation provides a powerful tool for representing quantum channels operating between nodes represented by density matrices (see Fig.~\ref{model}). This notation is especially useful for visualizing operations like trace distance and fidelity, which are related to the trace norm. We demonstrate this in Fig.~\ref{model1}.

\subsubsection{Assumptions} 
We assume the following rules regarding wavefunction eigenvectors in bra-ket notation 
\begin{align}
&|\psi\rangle\langle\psi| = |\psi\rangle \\
&\text {Prob} X(x) = \text{Tr} \left\{\Lambda_x \rho\right\} 
\end{align}
where $\text{Tr}$ represents Trace of a matrix and if $\Lambda_x\rho$, then we need to conduct a measurement of $x$ from $\rho$. We also assume that for every quantum channel we must have some classical measurement, such that,
\begin{align} \label{eq3}
    \rho \rightarrow \sum_x \operatorname{Tr}\left\{\Lambda_x \rho\right\} \underbrace{|x\rangle\langle x|}_{\substack{\text { classical } \\ \text { basis }}}
\end{align}
Assumption in (\ref{eq3}) is crucial for determining the upper/lower bounds on quantum channels when combined. For example, this is the case of super-activated channels, wherein we combine two or more zero capacity channels to generate nonzero capacity, and is applicable to some kinds of quantum nodes which act as network functional operators, like, quantum switches, repeaters etc.

\subsubsection{Formulation of Operational Form of TD} 
Next we start with the analytical representation of Trace Distance (TD) as,
\begin{align}\label{eq4}
    \|\rho-\sigma\|_1 = \sum|\lambda i|
\end{align}
where $\{\lambda i\}$ are the resultant eigenvalues of the mixed quantum state $\rho - \sigma$ over the connected quantum network. Next, we explore the operational meaning of TD. To use TD in the calculation of network capacity, we need to formulate the variational characterisation of TD as an optimization problem, the representation of which is provided as, $\max_{0 \leq \Lambda \leq\{\Lambda(\rho-\sigma)\}} \frac{1}{2}\|\rho-\sigma\|_1$. Here, we are optimizing over $\Lambda$ elements as long as the elements are less than those in $I$ but still positive, maximising over operators to compute a probability. Therefore, intuitively, we can say that ``probability difference" is analytically equivalent to TD.

Using the above thought, we manipulate (\ref{eq4}) to generate an inequality given by,
\begin{align}\label{eq5}
    \operatorname{Tr}\{\Lambda \rho\} \leq \|\rho-\sigma\|_2 + \operatorname{Tr}\{\Lambda \sigma\}
\end{align}
In (\ref{eq5}), if we pick a $\Lambda$ that is not an optimal measure, we can resort to approximation by computing probability. It is worth-mentioning here that TD is a useful metric of measure in quantum information theory (QIT) as it provides an operational interpretation of quantities which we need for ``Hypothesis testing" and interpret the inequality as either a game or a semi-definite program (SDP).

The inequality in (\ref{eq5}) can be solved using a variation of the CHSH (Clauser-Horne-Shimony-Holt) game \cite{} approach combined with elements of the Prisoner's Dilemma \cite{}, and a quantum measurement involving states $\rho$ and $\sigma$. The estimation errors in decision-making in the game can be defined in terms of  two types of errors, Type I and Type II, involving Bayesian probability and measurements. From the game perspective, we are dealing with quantum states shared between two parties, Bonnie and Clyde (B and C), and a detective who is trying to determine which states they are using based on measurements, where,
\begin{itemize}
    \item \emph{States $\rho$ and $\sigma$}: These are the quantum states shared between Bonnie and Clyde. The detective wants to determine which of these states is being used.
    \item \emph{Measurement $\Lambda$}: This is the measurement performed by the detective. The condition $\Lambda \rho + \Lambda \sigma = I$ indicates that the measurement operators $\Lambda \rho$ and $\Lambda \sigma$ are complementary and form a complete set of orthogonal projectors.
    \item \emph{Error Analysis}: involves guessing $\sigma$ using Bayesian probability. The parameter $r \in [0, 1]$ represents the probability of guessing $\rho$, and $1 - r$ represents the probability of guessing $\sigma$. The quantity $\operatorname{Tr}[U \sigma \rho]$ means one is calculating the trace overlap between the actual state $\sigma$ and the guessed state $\rho$ using the Unitary operator $U$.

\end{itemize}
The actual error rates will depend on the specific values of $r$, the properties of the states $\rho$ and $\sigma$, and the properties of the measurement operator $U$ to arrive at a simple formula for producing estimation of errors using methods of meoments on a random value of $r$, probability of which can be minimized using optimization.

\subsubsection{Error Analysis}
Maximum Likelihood Estimation (MLE) frames the measurement problem as an optimization task. We seek the parameters that best explain the measured eigenvalues (states $\rho$ and $\sigma$) while incorporating errors. Error regions are constructed around these optimal parameters using methods of moments. Mathematically, this can be represented as,
\begin{align}
    P_{\text {error }}\left(\frac{1}{2}, \rho, \sigma\right) =&~\text{Prob}\big[\text{error}[\rho]\text{Prob}[\rho] \nonumber\\
    &+ \text{Prob}[\text{error}[\sigma] \text{Prob}[\sigma]\big] \nonumber\\
    =&~\Tr[(I-U_\rho)] + \Tr[U_\sigma] \cdot(1-r)
\end{align}
with Trace Distance itself then expressed as;
\begin{align}
    TD_{\phi} = \frac{1}{2}\left(1+\frac{1}{2}\|\rho-\sigma\|_1\right) +/- P_{\text{error }}\left(\frac{1}{2}, \rho, \sigma\right)
\end{align}
The above is the equation that gives us a concrete way to calculate the trace distance. This is essential for practically measuring the difference between the initial and final states represented by density matrices $\rho$ and $\sigma$. Now, refer back to Fig.~\ref{model} where, we represent the derived concept so far in terms of tensor network diagram.

\section{Interpretation of an Operational Form of the Trace Distance}\label{S3}

\noindent We can interpret $P_{\text{error}}\left(\frac{1}{2}, \rho, \sigma\right)$ as follows: if $\rho, \sigma$ are indistinguishable, i.e. in a wave-function $\psi$ of probabilities, then the trace distance between them has to go to $0$ and $P_{\text{error}}=\frac{1}{2}$ (i.e. we have the measure of $\rho$ vs $\sigma$ as a perfect coin toss). On the other hand if $\rho$ and $\sigma$ are maximally distinguishable (i.e. they exist on orthogonal subspaces then T.D $=1$ and $P_{\text{error }} \rightarrow$ 0 (i.e. trending towards zero errors)

The operational trace distance formula provides a tool to quantify the degree of change in quantum systems undergoing non-unitary evolution (which involves incoherent processes like relaxation and dephasing).  This brings us to the concept of master equations. Master equations offer a theoretical framework to describe the time evolution of open quantum systems.  Fundamentally, they model how a system interacts with its environment, leading to processes like dissipation, decoherence, and relaxation. Mathematically, the master equation is a differential equation that describes the time derivative of the density matrix of the system. The density matrix, denoted by $\rho$, is a mathematical representation of the quantum state of the system. The master equation is typically written in the Lindblad form \cite{3spindecoherence}:
\begin{align}
\frac{{d\rho}}{{dt}} = -i[H, \rho] + \sum_k \left( L_k \rho L_k^\dagger - \frac{1}{2} \{L_k^\dagger L_k, \rho\} \right)
\end{align}
where, $H$ represents the system's Hamiltonian, which describes the energy of the system. The first term on the right-hand side, $-i[H, \rho]$, accounts for the unitary evolution of the system governed by the Hamiltonian.

The second term on the right-hand side captures how the system interacts with its environment, leading to effects like dissipation (loss of energy) or decoherence (loss of coherence).  This term includes a sum over ``collapse operators", denoted by $L_k$. These operators represent different channels through which the system loses coherence or energy to the environment.  In simpler terms, they act like quantum versions of noise sources we encounter in classical systems.The terms $L_k \rho L_k^\dagger$ and $L_k^\dagger L_k \rho$ quantify the dissipation and decoherence processes, respectively.

In the QuTiP Python library, the versatile $\mathsf{qutip.mesolve}$ function handles evolution under both the Schrödinger equation and master equations within quantum networks. This flexibility allows the user to choose the appropriate model. To analyze the evolution of trace distance, we focus on how density matrices change over time according to the Schrödinger equation using the QuTiP library. Let's assume that the Hamiltonian of the system is denoted by $H$. The time evolution of the density matrix $\rho$ is given by $\frac{d\rho}{dt} = -i[H, \rho]$, where $[,]$ denotes the commutator.

Now, let's consider two initial density matrices $\rho_1(0)$ and $\rho_2(0)$ at time $t=0$. The time evolution of these density matrices can be written as:
\begin{align}
    \rho_1(t) &= e^{-iHt} \rho_1(0) e^{iHt}\\
    \rho_2(t) &= e^{-iHt} \rho_2(0) e^{iHt}
\end{align}
Using these expressions, we can calculate the trace distance at time $t$ as:
\begin{align}
    D(\rho_1(t), \rho_2(t)) = \frac{1}{2} \text{Tr} |e^{-iHt} \rho_1(0) e^{iHt} - e^{-iHt} \rho_2(0) e^{iHt}|
\end{align}
To represent the Hamiltonian $H$ in terms of orthogonal states, we can decompose it as: $H = \sum_{n,m} E_{nm} |n\rangle\langle m|$, where $E_{nm}$ are the energy eigenvalues and $|n\rangle\langle m|$ are the corresponding projection operators. The states $|n\rangle$ and $|m\rangle$ are orthogonal, i.e., $\langle n|m\rangle = \delta_{nm}$, where $\delta_{nm}$ is the Kronecker delta.

\section{Numerical Results}\label{S4}

\begin{figure}
    \centering
\includegraphics[width=0.45\textwidth]{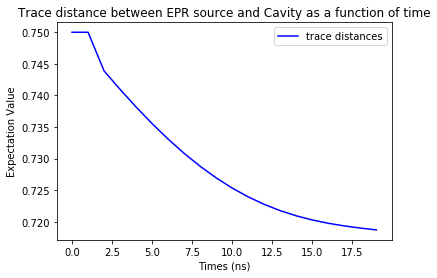}
%\centering
\vspace{-3mm}
\caption{Evolution of Trace distance measured between 2 nodes in a noisy environment as a function of time.}
\vspace{-5mm}
\label{fig2}
\end{figure}

\noindent The master equation uniquely empowers researchers to model both the unitary evolution of a density matrix and the impact of entanglement with environmental variables. This holistic approach elucidates the temporal dynamics of open quantum systems, including energy dissipation stemming from environmental interactions. The operators and rates within the master equation precisely characterize the nature and influence of these interactions. QuTiP's operational trace distance formula can be used within master equations to model noisy quantum environments. This allows us to track how a system evolves from a state of maximum difference (high initial trace distance) towards a less distinguishable steady state (trace distance approaching 1/2) due to decoherence.  This analysis provides insights into how distinguishable states become less so over time when coupled to a noisy environment (Fig.~\ref{fig2}).

The trace distance between an Einstein-Podolsky-Rosen (EPR) source and a cavity offers a way to quantify the distinguishability of quantum states within this system.  It could be used to measure how the ideal EPR state generated by the source deviates from reality due to imperfections, or to track the effects of cavity interactions on the entanglement of EPR photons.  Crucially, the specific choice of states being compared is important.  The complexity of modeling this system accurately might necessitate open quantum system approaches. Other relevant metrics include fidelity, which assesses similarity to a desired state, and entanglement measures, which directly characterize the entanglement present.
\begin{figure}
    \centering
\includegraphics[width=0.4\textwidth]{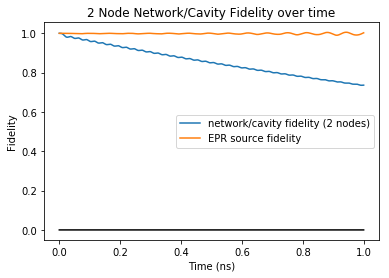}
%\centering
\vspace{-3mm}
\caption{Evolution of Link Fidelity between 2 nodes in a noisy environment as a function of time.}
\vspace{-5mm}
\label{fig3}
\end{figure}

We investigated the time-dependent fidelity of orthogonal states in a 2-node network in Fig.~\ref{fig3}. To establish a reference point, we introduced a 50/50 probability \emph{channel} between the network's \emph{supernode} (formed by the nearest and farthest nodes) and the EPR source.  Interestingly, the calculated fidelities of the supernode elements were very similar, suggesting a stronger connection within the supernode than between the supernode and the EPR source. Using Python's GRAPE (Gradient Ascent Pulse Engineering) algorithm, we modeled how the supernode's fidelity relative to the EPR source evolves over time. The model revealed noise-induced oscillations leading to periodic dephasing. The fidelity $\mathcal{F}$ between the target state $\rho_{\text{target}}$ and the evolved state $\rho_{\text{evolved}}$ can be calculated using the Quantum GRAPE technique in QuTip. For a SPDC (Spontaneous Parametric Down-Conversion) source, the fidelity can be expressed as:
\begin{equation}
\mathcal{F} = \frac{1}{2}(1 + \Tr\{\rho_{\text{target}} \rho_{\text{evolved}}\})
\end{equation}
The Frobenius norm can be used to relate the quantum nodes in tensor notation diagrams to the calculation of the quantum fidelity between two or more density matrices. Specifically, for two density matrices $\rho$ and $\sigma$, the quantum fidelity is given by: $F(\rho, \sigma) = \left(\text{Tr}\sqrt{\sqrt{\rho}\sigma\sqrt{\rho}}\right)^2$. The Frobenius norm can be used to simplify the calculation, and thus representation in tensor diagram, of the quantum fidelity. By representing the density matrices as tensors, the quantum fidelity can be expressed in terms of the Frobenius norm of the tensor difference between the two density matrices: $F(\rho, \sigma) = \left|\left|\rho^{1/2}\sigma\rho^{1/2}\right|\right|_F$, where $\rho^{1/2}$ represents the square root of the density matrix $\rho$. This relationship bridges the Frobenius norm (used in tensor analysis) with quantum fidelity. By calculating the Frobenius norm of the difference between two or more density matrix tensors, we directly assess their similarity.  This provides a quantifiable measure of quantum fidelity.
\begin{figure}
    \centering
\includegraphics[width=0.4\textwidth]{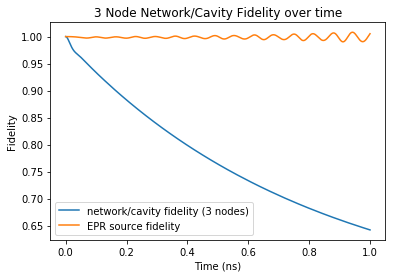}
%\centering
\vspace{-3mm}
\caption{Evolution of network Fidelity between 3 nodes in a noisy environment as a function of time.}
%\vspace{-5mm}
\label{fig4}
\end{figure}

Now we plot the evolution of 3-node network fidelity with time in Fig.~\ref{fig4}. While the overall network fidelity increases with time, the quality and coherence of the EPR source degrades over time. The reason can be attributed to noise and imperfections in the system resulting in degradation of the entanglement source, while still contributing to the improvement of network fidelity over time. In practice, the bridging together of fidelity and trace distance measures subject to errors produces a versatile benchmark, alluded to by previous work in quantum network benchmarking \cite{Helsen2023} but not developed with bounds to provide scope of the benchmark with errors in mind. We provide the bounds for the Trace distance vs Fidelity measures.

\subsubsection{Upper bound:}
The upper bound of the fidelity, denoted as $F_{\text{upper}}$, is given by: $F_{\text{upper}} = \sqrt{1 - F(\rho, \sigma)}$, where $F(\rho, \sigma)$ represents the fidelity between the states $\rho$ and $\sigma$. In terms of the Trace distance this is then: $F_{\text{upper}} = \sqrt{1 - \boldsymbol{TD}}$, where $\boldsymbol{TD}$ represents the trace distance between the two quantum states. The upper bound represents the maximum possible fidelity between the two states. It indicates the best-case scenario for the fidelity, assuming no errors or disturbances in the quantum network.

\subsubsection{Lower bound:}
The lower bound of the fidelity, denoted as $F_{\text{lower}}$, is given by: $F_{\text{lower}} = \sqrt{F(\rho, \sigma)}$, where $F(\rho, \sigma)$ represents the fidelity between the states $\rho$ and $\sigma$. In terms of the Trace distance this is then: $F_{\text{lower}} = \sqrt{\boldsymbol{TD}}$. The lower bound represents the minimum possible fidelity between the two states. It indicates the worst-case scenario for the fidelity, considering the maximum amount of errors or disturbances in the quantum system. The expected fidelity against trace distance bound measures are compared against the Fuchs-Van de Graaf benchmark [12] in Fig.~\ref{fig5}:

\begin{figure}
    \centering
\includegraphics[width=0.48\textwidth]{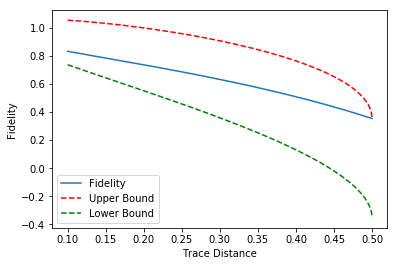}
%\centering
\vspace{-3mm}
\caption{Trace distance v/s fidelity bounds as compared to the Fuchs-Van de Graaf benchmark.}
\vspace{-5mm}
\label{fig5}
\end{figure}
From our perspective, the upper bound on fidelity serves as a gauge for the error level within a quantum channel.  A higher upper bound signifies a smaller overlap between the transmitted quantum states and the potential repeater's states. In simpler terms, this reduces the probability of information leaks that might occur due to errors or issues with synchronizing the measurement setup. The lower bound on fidelity, on the other hand, offers insights into potential communication system failures.  A very low lower bound suggests a situation where desynchronization has essentially caused the collapse of entanglement exchange. This collapse could be due to factors like depolarization within the channel (if we're using polarization density as the entanglement observable) or significant changes in the polarization extinction ratio (PER) between measured density matrices.  PER essentially measures the strength of a dominant polarization state compared to its orthogonal counterpart after traveling through a system.

\section{Conclusion}\label{S2}

Our operational form of trace distance, calculated between two nodes within a quantum network, provides a valuable metric for assessing overall network fidelity as the network scales. By introducing error bounds on fidelity vs. trace distance, we establish a robust benchmark.  We propose grounding this benchmark in  physical observables, such as polarization density (represented on the Poincaré sphere) and polarization extinction ratio (PER). This approach is applicable to experimental setups using polarization-entangled light in both fiber and free-space networks, with links to established performance metrics like the Quantum Bit Error Rate (QBER).

\end{document}